\documentclass[12pt.sty]{article}
\usepackage{amsfonts, amssymb, amsmath, amsthm, }
\theoremstyle{plain}
    \newtheorem{Th}{Theorem}[section]
    \newtheorem{Lem}[Th]{Lemma}
    \newtheorem{Prop}[Th]{Proposition}
    \newtheorem{Cor}[Th]{Corollary}
    \theoremstyle{definition}
    \newtheorem{Ex}[Th]{Example}
    \newtheorem{Def}[Th]{Definition}
    \newtheorem{Rem}[Th]{Remark}
\begin{document}
\title{Multiplicative Aspects of the Halperin-Carlsson Conjecture}

\author{Volker Puppe}

\date{}

\maketitle


\begin{abstract}
\noindent We use the multiplicative structure of the Koszul
resolution to give  short and simple proofs of some known
estimates for the total dimension of the  cohomology  of spaces
which admit free torus actions and some analogous results for
filtered differential modules over polynomial rings. We also point
out the possibility of improving these results in the presence of
a multiplicative structure on the so-called minimal Hirsch-Brown
model for the equivariant cohomology of the space.
\end{abstract}



\section{Introduction}


All spaces considered in this note are assumed to be paracompact
Hausdorff.\\

\noindent The Halperin-Carlsson conjecture is the following
statement:

\begin{quote}
$(\bf{HCC})$ \hspace {1mm}If a torus $S^1 \times  \dots \times
S^1$ (resp. a $p$-torus, $\mathbb Z/p  \times  \dots  \times
\mathbb Z/p, p$ prime) of rank $r$ can act freely on a finite
dimensional space $X$, then the total dimension of its cohomology,
$ \sum_i dim_k H^i(X;k) \geq 2^r$. Here $k$ is a field
characteristic $0$ (resp. $p$).
\end{quote}

Many authors, see e.g.
[Al],[AH],[AP1],[AP2],[AB],[Ba],[Bo],[Ca1]-[Ca5],[Hl],[Hn],[Ya],
have studied variants of this problem and have contributed results
with respect to different aspects of the conjecture. For part of
the literature and a summary of certain results see e.g. the
discussion in [AD] and [AP2],(4.6.43).

It seems
that in those cases where one can prove the conjecture, the
multiplicative structure plays an essential role, or one considers
rather special spaces like
products of spheres of the same dimension. \\
To my knowledge the best general estimate for $ \sum_i dim_k
H^i(X;k)$ is the following. We refer to the real torus case as
'case(0)', and to the $p$-torus as 'case(p)'.

\begin{Th}
(a) In case(0) the conjecture holds for $r \leq 3$, and for any
$r\geq 3$ one has: $ \sum_i dim_k H^i(X;k) \geq
2(r+1)$.\\
(b) In case(2) the conjecture holds for $r \leq 3$, and for any
$r\geq 3$: \hspace{0mm} $ \sum_i dim_k H^i(X;k) \geq
2r$.\\
(c) In case(p), $p$ odd, the conjecture holds for $r \leq 2$, and
for any $r$:  $ \sum_i dim_k H^i(X;k) \geq
r+1$.\\
\end{Th}

In [ABI] and [ABIM] related results in a rather general
algebraic context are discussed and proved.\\

In this note me make use of some rudiments of multiplicative
structure, which come from the $dga$-structure of Koszul
complexes, to give simple proofs of some of the known results
mainly in the graded context. The method of proof suggests on one
hand that it might be useful to study the multiplicative structure
up to homotopy on the minimal Hirsch-Brown in more detail; on the
other hand the method seems flexible enough to apply also in more
general
algebraic situations (cf. Cor.4.6).\\

\section{Maps between Koszul complexes}

In this section $k$ denotes a field of arbitrary characteristic.
For any integer $m$ let $K_r(m)$ be the Koszul complex
corresponding to the regular sequence $(t_1^{m+1},...,t_r^{m+1})$
in $R:= k[t_1,...,t_r]$, ($deg(t_i) = 1$). $K_r(m)=\Lambda _R
(s_1^m,...,s_r^m)$, is the exterior algebra over $R$ in $r$
generators $s_1^m,...,s_r^m$ of degree $m$ with differential
$d({s_i}^m) = t_i^{m+1}, d(t_i) = 0, i = 1,...,r$, extended as a
derivation to obtain a $dga$(differential graded algebra). The
Koszul complex can be considered as the minimal resolution of the
$R$-module $R/(t_1^{m+1},...,t_r^{m+1})$. We consider $R$-linear
maps $\gamma: K_r(m) \longrightarrow K_r(0)$ from the $R$-cochain
complex $K_r(m)$ to the $R$-cochain complex $K_r(0)$, which lift
the projection $ H(K_r(m)) = R/(t_1^{m+1},...,t_r^{m+1})
\longrightarrow
 H(K_r(0)) = R/(t_1^{1},...,t_r^{1}) \cong k.$ But we do not assume that the
maps preserve the length of exterior products nor the grading,
i.e. they are just morphisms of the underlying differential
modules. By $rk(\gamma)$ we mean the rank of the map induced by
$\gamma$ on the localized modules, inverting all non-zero
homogenous polynomials in $R$.

\begin{Lem}

(a) If $\gamma$ is a $dga$-map (in particular multiplicative),
then $rk(\gamma) = 2^r$. (cf. [BE], Prop.1.4.) \\
(b) For any $\gamma$ as above $rk(\gamma) \geq 2r $.\\
\end{Lem}

\noindent {\bf Proof:}  Since $K_R(0)$ is a free resolution of
$k$, considered as a $R$-module via the canonical augmentation,
which maps all the $t_i$'s to $0$, any two maps of the form
$\gamma$ are homotopic. In particular any such $\gamma$ is
homotopic to the map $\iota$ defined as the dga-map which sends
$s_i^m$ to $t_i^m s_i^0$. Hence $\gamma$ induces the same map as
$\iota$ in cohomology for any coefficients. We consider the
induced map in cohomology with coefficients in $ \bar{R} :=
R/(t_1^{m+1},...,t_r^{m+1})$ and use the map $\iota$ for
calculations. The induced boundary on $K_r(m) \otimes \bar{R}$
vanishes and it is easy to see, that $H(K_R(0) \otimes \bar{R})$
is the exterior algebra over $k$, generated by $t_1^m
s_1^0,...,t_r^m s_r^0$ (Note that the superscript for $s_i$ is
used just to identify the element, whereas the superscript for
$t_i$ denotes an exponent.) The induced map ${\gamma}^* =
{\iota}^*: K_r(m) \otimes \bar{R} \longrightarrow   H(K_r(0)
\otimes \bar{R})$ maps $t_i$ to $[t_i^m s_i^0]$ and is
multiplicative. In particular the element $s^m_{1...r} := s_1^m
\wedge ...\wedge s_r^m$ is mapped to the non-zero element $[t_1^m
\cdots t_r^m s^0_{1...r}]$. Therefore $\gamma (s^m_{1...r})$ is
non-zero in $K_R(0)$. If $\gamma$ is assumed to be multiplicative
then it follows that it must be injective for the following
reason: For any element $x \in K_r(m)$ there exists an element $y
\in K_r(m)$ such that $xy = qs^m_{1...r}$, with $ q \in R$
(because the exterior algebra fulfills Poincar\'e duality over the
quotient field). So for $x \neq 0$, $\gamma (x)$ must be non-zero,
otherwise $\gamma (xy)$ would vanish, which is not the case since
$\gamma (s^m_{1...r})$ is non-zero in $K_R(0)$.  Hence we get
part(a) of the lemma, since localization is exact.\\
To prove part(b) we change $\gamma$ by first restricting the map
to the free $R$-submodule $\Lambda^{\leq 1}$ of $K_r(m)$ generated
by $ 1, s_1^m,...,s_r^m $ and then extending it multiplicatively.
On this submodule the two maps coincide and are injective by
part(a). Since $\gamma$ commutes with the respective boundaries it
is also non-zero on all elements in $K_r(m)$, which are mapped by
the boundary into non-zero elements of $\Lambda^{\leq 1}$. An
$R$-linear combination of the elements $
1,s_1^m,...,s_r^m,s_{12}^m,...,s_{2r}^m$ is mapped to a non-zero
element in $\Lambda^{\leq 1}$ if at least on of the coefficients
of the $s_{1j}^m$ is non-zero. Hence $\gamma$ is injective on the
free $R$-module generated by th $2r$ above elements.
This gives part(b) of the Lemma. $\hfill\Box$\\

Some time ago Martin Fuchs found an example of a map $\gamma$ with
$rk(\gamma) < 2^r$, and $r = 4$. The following example shows that
one can not improve the estimate in part(b) of the above lemma for
$r = 3$ without additional assumptions.

\begin{Ex}
The following map $\gamma : K_3(1) \longrightarrow K_3(0)$ is
homotopic to the standard map $\iota$ and has rank equal to
$2r=6$, ($k = \mathbb F_2$). Define $\gamma = \iota + dh+hd$ with
$h(s_1^1) = s^0_{123},  h(s^1_{23}) = t_3 s^0_{12}$ and otherwise
equal to $0$ for the standard basis of $K_3(1)$. Direct
computation shows that $\gamma$ vanishes on $x:= t_1t_3s^1_{12} +
t_2s^1_{123}$ and on $dx$, and $x$ and $dx$ are linearly
independent over $R$. Hence $rk(\gamma) \leq 6$ and by the above
lemma it must be indeed equal to $6$. $\hfill\Box$\\
\end{Ex}
In view of the topological applications, (see Section 5), we will
also consider the degree convention $deg(t_i) = 2$ in case $(0)$;
Lemma 2.1. holds in this context, too. \\
The above example
(with appropriate signs) also works for other fields if $deg(t_i)
= 1$, but it does not work, if one
puts $deg(t_i) = 2$ (cf. Theorem 5.1.(b)).  \\

For later use we prove the following technical result:

\begin{Lem}
In case (0), (i.e. $deg(t_i) = 2$, $deg(s_i^m) = (2m+1)$), if $r
\geq 3$ then any map $\gamma$ as above, which preserves degrees,
is injective on the free $R$-submodule of $K_r(m)$, which is
generated by the elements {$s_1^m, ..., s_r^m, s^m_{123}$}.
\end{Lem}

\noindent {\bf Proof:} Let $x = qs^m_{123} + \sum_1^3 q_i  s_i^m$
with $q, q_i \in R$ be a homogeneous element in $K_r(m)$. Note
that we already know from the above arguments, that we only need
to consider the case $q \neq 0$. Let $J$  be the ideal generated
by $qt_i^{m+1}, i = 1,...,r$ and $q_i t_i^{m+1}, i = 1,2,3$. We
consider cohomology with coefficients in $R/ J$, so $x$ becomes a
cycle; and we show that $\iota ^*$ and hence $\gamma ^*$ is
non-zero on $[x]$, which, of course implies that $\gamma (x) \neq
0$. The maps $\iota$ sends the element $x$ to $qt_1^m t_2^m t_3^m
s^0_{123} + \sum_1^3 q_i t_i^m s_i^0$. Since $deg(s_i^m) = 2m+1$
in case (0), one has $deg(q_i) = deg(q)+4m+2$. So $deg(qt_1^m
t_2^m t_3^m) < deg(q_i t_i^{m+1}$, and since $qt_1^m t_2^m t_3^m$
is not contained in the ideal generated by $qt_i^{m+1}, i =
1,...,r$, it is also not contained in
 $J$.  Hence $[\iota (x)]$ is non-zero in $H^*(K_r(0)
\otimes_R  R/ J)$.         $\hfill\Box$\\

Lemma 2.3. also holds (for any characteristic) if $deg(t_i) = 1$,
but we will apply it later in the above form.

\section{Minimal models}

We give a brief account of some facts about (additive) minimal
models of cochain complexes over the graded polynomial ring
$R=k[t_1,...t_r] $, where $k$ is a field of arbitrary
characteristic. The (additive) minimal model here plays a role
similar to that of the minimal resolution in homological algebra.
We do not consider any product structure (besides the $R$-module
structure) on the complexes in this section. So these minimal
models should not be confused with the (multiplicative) Sullivan
minimal models. The material is known, see e.g.[AP2], Appendix B,
and can also be drawn from several other sources which sometimes
deal with much more general situations, but we hope that it will
be convenient for the reader to get a short and rather elementary
presentation of what is needed in this note. \\

Let $\tilde C$ be a free cochain complex over $R$ with boundary
$\tilde d$ of total degree 1.  As $R$-module, $\tilde C \cong C
\otimes R$, where $C$ is a $k$-vector space, and the unspecified
tensor product,$\otimes$, is taken over the field $k$. We want to
define a free cochain complex, which is homotopy equivalent to the
given one and has minimal rank over $R$. Such a complex turns out
to be unique up to isomorphism and is called the (additive)
minimal
model of $\tilde C$.\\

The $R$-linear boundary $\tilde d$ can be written as a sum $d+d'$,
with $d$ a boundary on $C \otimes k$, and $d'(\tilde C) \subset C
\otimes  I$, where $I$ denotes the augmentation ideal of R. If we
consider the elements in $\tilde C$ as polynomials in the
variables $t_i, i= 1,...,r$ with coefficients in $C$, the part $d$
of the boundary correspond to the coefficients of $1$ i.e. the
constant part of $\tilde d$, and $d'$ to the higher
terms in the $t_i$'s.\\
$(C,d) \cong (\tilde C \otimes_R k, \tilde d \otimes_R k)$ , where
$k$ is an $R$-module via the augmentation, is a cochain complex
over the field $k$. Hence we can write $C$ as a (non-canonical)
direct sum $H \oplus B \oplus D$, where $H = H(C,d)$, and $d$
corresponds to an isomorphism (also called $d$) from $D$ to $B$.
We extend this isomorphism $R$-linearly to $B \otimes R$ and
obtain this way a contractible cochain complex $\tilde N \cong (B
\oplus D) \otimes R$. We imbed this complex into $C \otimes R$,
imbedding $B \otimes R$ by $R$-linear extension of the embedding
of $D$ into the direct sum above and extending this map to be
compatible with the respective boundaries. Note that this is well
defined since the boundary in $\tilde N$ maps $D \oplus R$
isomorphically to $B \otimes R$, but it is not(!) the $R$-linear
extension of the embedding of $B \oplus D$, if $d'$ is non-zero.
Nevertheless the map defined is an embedding of a contractible
cochain complexes onto a direct summand as an $R$-module, and the
quotient complex, $ \tilde H$, is isomorphic as an $R$-module to
$H \otimes R$, but inherits  a twisted boundary $\bar{d}$. If one
tensors $\tilde C \cong \tilde H \oplus \tilde N$ with $k$,
considered as a $R$-module via the augmentation (or in other
words: if one restricts to constant terms with respect to the
variables $t_i$), then one gets back the direct sum decomposition
given above. In particular the constant part of the boundary in $
\tilde H$ is zero.\ The complex $ \tilde H$ is $R$-homotopy
equivalent to $ \tilde C$, since it is obtained from the latter by
dividing out a contractible direct summand. It is not difficult to
see that up to isomorphism of $R$-complexes there is only one free
$R$-complex which is $R$-homotopy equivalent to $ \tilde C$, and
has a boundary with vanishing constant term (i.e. which vanishes
when tensored with $k$ over $R$). The complex $ \tilde H$ is
called the (additive) minimal model of the complex  $ \tilde C$.
It follows from ${\tilde d} \circ {\tilde d} = 0$ that the part of
${\tilde d}$, which is linear in the $t_i$'s, anti-commutes with
$d$. Therefore this linear part induces  a map on $H(C,d)$ which
is in fact the linear part of the boundary of $\tilde H$. For
every index $i$ we define a map ${\lambda}_i$ from $H$ to itself
by assigning to an element $x \in H$ the coefficient of $t_i$ in
$\bar{d}(x)$. Since ${\bar{d}} \circ {\bar{d}} = 0$, the maps
${\lambda}_i$ anti-commute for $i,j = 1,...,r$. Hence they define
an action of the exterior algebra ${\Lambda}_k
({\lambda}_1,...,{\lambda}_r)$ on $H$. Let $\Lambda ^+$ be the
ideal generated by ${\lambda}_1,...,{\lambda}_r$. The length,
$\ell _{\Lambda}(H^q)$, of $H^q$ as a ${\Lambda}_k
({\lambda}_1,...,{\lambda}_r)$ - module is definite as the minimal
integer $i$, such that $(\Lambda ^+)^i H^q = 0$. One has proper
inclusions $(\Lambda ^+)^i H ^q\supset (\Lambda ^+)^{i+1} H^q $
for $i= 0,...,(\ell _{\Lambda}(H^q)-1)$.
\\
Assume now that $dim_k H$ is non-zero and finite. We define a
filtration on the minimal model $\tilde H$ by subcomplexes
inductively in the following way (cf.[AP2],Sec.1.4):

\begin{Def}
Let $\tilde H \cong H\otimes R$ be a minimal model as above. We
define\\ ${\cal F} _0(H) = 0$ and ${\cal F} _0(\tilde H)= 0$,\\
${\cal F} _1(H): = ker ({\bar d}:H \cong {H \otimes k}
\longrightarrow \tilde H$) and ${\cal F}_1(\tilde H): = {\cal
F}_0(H) \otimes R$. \\Let us assume that ${\cal F}_i(\tilde H): =
{\cal F}_i(H) \otimes R$ is already defined, then we put\\ $ {\cal
F} _{i+1}(H) : = {\bar d}^{-1} ( {\cal F}_i (\tilde H)) \cap (H
\otimes k) $ and ${\cal F}_{i+1}(\tilde H): = {\cal F}_{i+1}(H)
\otimes R$.\\ The length of the filtration, $\ell ({\cal F} _*
(\tilde H))$, is the smallest index $i$ for which ${\cal F}_i
(\tilde H)$ and ${\cal F}_{i+1} (\tilde H)$ coincide.
\end{Def}

We summarize some of the properties of this filtration in the
following proposition.

\begin{Prop}
(a) One has proper inclusions of free subcomplexes \\
${\cal F} _0(\tilde H) = 0 \subset {\cal F} _1(\tilde H)
\dots\subset {\cal F} _{\ell ({\cal F} _* (\tilde H))} (\tilde H)
= \tilde H$,
which are direct summands as $R$-modules.\\
(b) $\bar d ({\cal F} _i (\tilde H)) \subseteq {\cal F} _{i-1}
(\tilde H)$, for $i = 0,...,\ell ({\cal F} _* (\tilde H))$\\
(c) The complex $\tilde H$ admits an augmentation
$\varepsilon:\tilde H \longrightarrow k$ compatible with the
respective boundaries (where the boundary on $k$ is trivial), such
the restriction to ${\cal F} _0 (\tilde H)$ is surjective.\\
(d) The boundary $\bar d: {\cal F} _i (\tilde H) \longrightarrow
{\cal F} _{i-1} (\tilde H)$ induces a non-trivial map from $
({\cal F} _i (\tilde H))/ {\cal F} _{i-1} (\tilde H)$ to $({\cal
F} _{i-1}
(\tilde H))/ {\cal F} _{i-2} (\tilde H)$ for $i = 1,...,\ell ({\cal F} _* (\tilde H))$.\\
(e) $dim_k H \geq\sum_q \ell _{\Lambda}(H^q) \geq \ell ({\cal F}
_*
(\tilde H))$\\
(f) If the action of $\Lambda$ on $H$ is trivial, i.e. the terms
in $\bar d$, which are linear with respect to the $t_i$'s vanish,
then $\sharp$\{$q; H^q \neq  0$\}  $\geq \ell ({\cal F} _* (\tilde H))$.\\
\end{Prop}

\noindent {\bf Proof:} The properties (a)-(e) follow directly from
the definition using in particular the following facts: \\
- the total degree of the boundary is 1,\\
- the constant part vanishes,\\
- the linear part of the boundary of $(\Lambda ^+)^i H^q \otimes
k$ is contained in $(\Lambda ^+)^{(i+1)} H^q \otimes R$,\\
and the higher order parts are contained in $H^{<q} \otimes R$. \\
Part (f) follows immediately from the fact that, under the
assumption made, the boundary of $H^i \otimes k$ is contained in
$H^{\leq {(i-1)}} \otimes R$ for degree reasons. $\hfill\Box$\\

The boundary of $\tilde H$ induces a boundary on the associated,
graded complex of the filtration $\cal F _*$, which is non-trivial
for all $({\cal F}_i (\tilde H)/{\cal F}_{i-1} (\tilde H)), i =
0,... \ell (\cal F _*)$. Since the composition of two successive
boundary maps vanishes, one has $rk({\cal F}_i (\tilde H)/{\cal
F}_{i-1} (\tilde H)) \geq 2$ for $i = 1,..., (\ell ({\cal F}
_*)-1)$, otherwise one of the two maps in the composition would
have to vanish, which is not the case, see Prop.3.2.(d). Hence we
get(cf.[AP2],Cor.(1.4.21):

\begin{Cor}
$dim_k(H^*)\geq  2(\ell ({\cal F} _* (\tilde H))-1)$.
\end{Cor}

\section{Factorization}

In this section we combine the results of two previous sections.
We consider free cochain complexes $\tilde C$ over the ring $R$ as
in the previous section. We assume that $\tilde C$ has an
augmentation $\varepsilon: \tilde C \longrightarrow k$, which
induces a surjection in cohomology.

\begin{Prop}
Let $\tilde C$ be a cochain complex over $R$ as above and such
that $H^{\geq {(m+1)}}(\tilde C )$ vanishes, then there exists a
map of $R$-complexes $\alpha:K_r(m) \longrightarrow \tilde C$,
such that $\varepsilon ^* \alpha ^*:H(K_r(m)) =
R/(t_1^{m+1},...,t_r^{m+1}) \longrightarrow k$ is the canonical
projection.
\end{Prop}

The proof of this proposition is by standard homological algebra.
One defines $\alpha$ inductively over the lengths of exterior
products in $K_r(m)$. The assumption on the vanishing of the
cohomology in high degrees allows to choose the images of the
elements $s^m_i$ compatible with the boundary, and so on. Cf. e.g.
[AP2],Lemma (1.4.17). $\hfill\Box$\\

Let $\tilde C$ be free differential $R$-module with a filtration
${\cal F}_*(\tilde C)$, which has the properties (a) and (b) of
Prop.3.2.(for $\tilde C$ in place of $\tilde H$). In [ABI] such an
object is called a free differential flag. We assume in addition
that $\tilde C$ has an surjective augmentation $\varepsilon: {\cal
F}_0 (\tilde C) \longrightarrow k$, which extends to morphism of
differential modules $\varepsilon:\tilde C \longrightarrow k.$ We
consider $K_r(0)$ with the filtration by length of exterior
products and with the canonical augmentation.

\begin{Prop}
Under the above assumptions there exists a morphism of
differential $R$-modules $\beta: \tilde C \longrightarrow K_r(0)$,
which commutes with the respective augmentations.
\end{Prop}

Again the proof is by standard homological algebra using induction
over the filtration degree. See [Ba], cf.
[AP2], Lemma (1.4.18)),  for the corresponding result in the presence
of an additional grading. $\hfill\Box$\\

\begin{Cor}
If $\tilde C$ is a free cochain complex over the graded polynomial
ring $k[t_1,...,t_r]$, $k$ a field of arbitrary characteristic,
$deg(t_i) = 1$, such that $dim_k H^*(\tilde C)$ is non-zero finite, then\\
(a) $dim_k H = \sum _i H^i(\tilde C \otimes k)  \geq 2r$,\\
(b) For the length of the filtration on the minimal model of
$\tilde C$ one has:\\
$\ell ({\cal F} ^*( \tilde H)) \geq \sum_q \ell_{\Lambda} (H^q) \geq(r+1)$.\\
(c) If the linear part of the boundary on the minimal model
vanishes, one has: $H^q(\tilde C)$ is non zero for at least $r+1$
degrees $q$.
\end{Cor}

\noindent {\bf Proof:} Part(a): By the two above propositions,
applied to the minimal model $\tilde H$ of $\tilde C$, one obtains
morphisms $\alpha$ and $\beta$ such that the composition $\beta
\alpha$ sends $1 \in K_r(m)$ to $1 \in K_r(0)$. So, by the Lemma
2.1, the rank of this composition is greater or equal to $2r$.
Since the map factors through the minimal model, $\tilde H$, the
rank of the model (when localized) must also be $\geq 2r$. But as
an $R$-module this model is free of dimension $ dim_k H = \sum _q
H^q(\tilde C \otimes k)$. Hence the assertion follows.\\
Part(b): The length of the filtration by exterior products on
$K_r(0)$ is equal to $(r+1)$. The element $s^0_{1...r}$ has
filtration length $r$. The  multiple $t_1^m\cdots t_r^m
s^0_{1...r}$ represents a non-zero element in $H^*(K_r(0) \otimes
_R R/J)$ . This class is in the image of $\iota ^*$ and hence also
in the image of $\beta ^*$ (see Section 2). Since by Prop. 4.2.
the map  $\beta$  can be chosen to preserve filtrations, it
follows that the length of any filtration of $\tilde H$ has
to be strictly greater than r. Together with Prop.3.2.(e) one obtains the assertion.\\
Part(c) follows from Prop.3.2.(f).   $\hfill\Box$\\

\begin{Rem}
Instead of using Lemma 2.1 one can deduce part(a) of the above
corollary from part(b) and Cor.3.3.
\end{Rem}

\begin{Rem}
There are far reaching recent generalizations, which give similar
bounds in a much more general algebraic context (see [ABI],
[ABIM]. Our main point here is to present a rather elementary
proof in the most classical, graded situation. But the method has
also some potential to be applied more generally, see Cor.4.6.  On
the other hand it is rather doubtful that it can lead to better
estimates without substantial additional effort, as the above
Example 2.2. shows.
\end{Rem}

Let $(\tilde C ,\tilde d)$ be an (ungraded) free differential
$R$-module ($R = k[t_1,...,t_r]$ also considered without grading)
with a filtration and an augmentation as above. We assume that the
filtration is minimal, i.e. that also Prop.3.2.(d)(for $\tilde C$
in place of $\tilde H$) holds. Finally we suppose that the
annihilator ideal of $H^*(\tilde C)$ contains the elements
$t_1^{m+1},...,t_s^{m+1}$. This is a replacement for the
assumption that $H^{\geq {(m+1)}}(\tilde C )$ vanishes in the
graded case.

\begin{Cor}
Under the above assumptions one has:\\
(a) $rk_R (\tilde C) \geq 2s$.\\
(b) The length of the filtration of $\tilde C$ is at least $s+1$.
\end{Cor}

The proof is completely analogous to that of Cor.4.3.(a) and (b).
The assumptions above allow to obtain a factorization up to
homotopy of the standard map $\iota:K_s(m) \longrightarrow K_s(0)$
through the filtered differential module $\tilde C$ using the
Propositions 4.1. and 4.2., more precisely: The assumption on the
annihilator ideal allows to apply Proposition 4.1. adapted to the
situation at hand, and the assumptions on the augmentation make
sure that one can apply Proposition 4.2. (We shift from
$k[t_1,...t_r]$-modules to
$k[t_1,...,t_s]$-modules via the canonical inclusion and projection.)$\hfill\Box$\\

The following corollary is of a somewhat different nature.\

\begin{Cor}
If $\tilde C$ is as in Cor.4.3, then\\
$dim_k H({\tilde C} \otimes_R  R/(t_1^{m+1},...,t_r^{m+1}))
\otimes_R R/(t_1,...t_r)) \geq 2^r$, for large enough m. In other
words: The minimal number of generators of $H({\tilde C} \otimes_R
R/(t_1^{m+1},...t_r^{m+1}))$ as an $R$-module is greater or equal
to $2^r$.
\end{Cor}

\noindent {\bf Proof:} The map ${\iota}^*:H^*(K_r(m) \otimes_R
R/(t_1^{m+1},...,t_r^{m+1})) \longrightarrow   H^*(K_r(0)
\otimes_R R/(t_1^{m+1},...,t_r^{m+1}))$ maps $[t_i]$ to $[t_i^m
s_i^0]$, see proof of Lemma 2.1. Taking the tensor product of this
map with $R/(t_1,...,t_r)$ over $R$ one gets an isomorphism. Since
the map factors through $H({\tilde C} \otimes_R
R/(t_1^{m+1},...t_r^{m+1}))
\otimes_R R/(t_1,...t_r))$, the assertion follows. $\hfill\Box$\\

Note that the above corollary also applies to the minimal model
$\tilde H$ of $\tilde C$. So \\
$dim_k H({\tilde H} \otimes_R R/(t_1^{m+1},...t_r^{m+1)}))
\otimes_R R/(t_1,...t_r)) \geq 2^r$. But the Halperin-Carlsson
conjecture in this algebraic context can be stated as \\
$dim_k (H^*(\tilde H \otimes_R R/(t_1,...,t_r)) = \tilde H
\otimes_R R/(t_1,...,t_r) \geq 2^r$. Although these two statements
look rather similar, they differ by an interchange of taking
homology and tensor product.

\section{Applications to torus actions}

In this section we give applications of our previous results to
free torus (resp. 2-torus) actions. We then have $H^*(BG;k) \cong
R = k[t_1,...,t_r]$, where $char(k) = 0$ and $deg(t_i = 2)$ (resp.
$char(k) =2$ and $deg(t_i) = 1$). \ Let $X$ be a finite
dimensional space on which a torus, $S^1 \times...\times S^1$
(resp. a 2-torus, $\mathbb Z/2 \times ...\times \mathbb Z/2 $) of
rank $r$ acts freely.\ We will use Cor.4.3. to get estimates for
the size of the cohomology of $X$ with coefficient in a field $k$
of characteristic zero (resp. 2). We consider the equivariant
cohomology of the space $X$ and briefly recall some facts about
the minimal Hirsch-Brown model of the equivariant cohomology for a
$G$-spase $X$ (see [AP2] for details).  For a $G$-space $X$ the
Borel construction gives a fibration $X \longrightarrow X_G : = EG
\times _G X \longrightarrow BG$ where $BG$ is the classifying
space and $EG$ the universal (free, contractible) $G$-space. For
$G=(S^1)^r$ (resp. $(\mathbb Z/2)^r$) there is the following
additive minimal cochain model over $R \cong H^*(BG;k)$ for the
cohomology of $X_G$, the so-called minimal Hirsch-Brown model
(s.[AP2]): $H^*(X;K) \tilde{\otimes} H^*(BG;k)$, where the tilde
indicates that the tensor product carries a twisted differential,
which in a sense
reflects the $G$-action on the cochain level (cf.[AP2]).\\

In case $G = (\mathbb Z/p)^r$ the minimal Hirsch-Brown model is in
general not a cocchain complex over $H^*(BG;k)$, but only a module
over the polynomial part of $H^*(BG;k)$. The behavior with respect
to the exterior part of $ H^*(BG;k)$ is rather complicated and our
algebraic results above do not suffice to give the results stated
in the introduction for p-tori, $p$ odd. We refer to [Ba] and
[AP2],(1.4.14) for more involved proofs of these results in a
similar spirit.\\

If $G$ acts freely on $X$ then $X_G  \simeq  X/G$; in particular\
\begin{center}
$H^*(X_G;k) = H^*(H^*(X;k) \tilde{\otimes} H^*(BG;k))$ $ \cong
H^*(X/G;k)$
\end{center}
In the case $G =
(S^1)^r$ the last isomorphism even holds under the weaker
assumption that all isotropy groups are finite. If, in addition,
$X$ is finite dimensional then so is $X/G$, and hence $H^i(X_G;k)
\cong H^i(X/G;k) = 0$ for large enough i, say $i \geq 2(m+1)$ in
case $G = (S^1)^r$ (resp. $i \geq
m+1$ in case $G = (\mathbb Z/2)^r$).\\

We apply Cor.4.3 and Lemma 2.3 to the minimal Hirsch-Brown model
to obtain part of the result stated in the introduction.

\begin{Th}
(a) If $X$ is a finite dimensional space on which a 2-torus $G =
(\mathbb Z/2)^r$ acts freely, then $dim_k \sum_i H^*(X;k) \geq
2r$.\\
(b) $\sum_q \ell_{\Lambda} (H^q(X;k)) \geq r+1$; in particular, if
the action on $X$ induces the trivial action in cohomology, then
$H^i(X;k)$ is non-zero for at least $(r+1)$ degrees $i$ ($char(k)
= 2$). \\
(c) If $X$ is a finite dimensional space on which a torus $G =
(S^1)^r$ acts almost freely, then $dim_k \sum_i H^*(X;k) \geq 2r$.
The number of degrees, for which $H^i(X;k)$ is non-zero is at
least $r+1$. \\
For $r \geq 3$,  $dim_k \sum_i H^*(X;k) \geq 2(r+1)$ ($char(k) =
0$).
\end{Th}

\noindent {\bf Proof:} Part(a) and the first two parts of (b)
follow immediately from Cor.4.3. (Note that trivial action in
cohomology implies that the linear part of the boundary of the
minimal Hirsch-Brown  model  vanishes in case(2); for case(0) this
is always true). To show the slightly improved inequality in
part(b) we observe that, in case(0), $\tilde H$ when localized (by
inverting all non-zero homogeneous polynomial in $R$) inherits a
$\mathbb Z/2\mathbb Z$-grading by even and odd degree. Since the
localized cohomology vanishes, the ranks of the even and the odd
part in $\tilde H$ are equal. Now Lemma 2.3 together with the
above factorization shows that the odd part of $\tilde H$ has rank
at least r+1. So  $rk (\tilde H) = dim_k \sum_i H^i(X;k) \geq
2(r+1)$.                                      $\hfill\Box$\\

 The minimal Hirsch-Brown model carries a multiplicative
structure which induces the cup product in equivariant cohomology.
In general the multiplication on the model is commutative and
associative only up to (higher) homotopies. In rather special
cases, e.g. in case(0) when  X  is a product of spheres of odd
dimension, the minimal Hirsch-Brwon model coincides with the
Sullivan minimal model of the Borel construction $X_G$, and hence
is a differential graded algebra (dga). In such a situation one
can derive from the above results the following corollary.

\begin{Cor}
If the minimal Hirsch-Brown model of the finite dimensional, free
$G$ -space $X$ carries a dga-structure (over $R$), then \ $dim_k
\sum_i H^*(X;k) \geq 2^r$ (cf.[BE], Prop.1.4)
\end{Cor}
\noindent {\bf Proof:} Under the given hypothesis the map $\alpha
: K_r(m) \longrightarrow \tilde H $ can be chosen to
multiplicative, and hence an argument analogous to that given for
Lemma 2.1 (a) shows that $\alpha$ must be injective. Therefore
$rk(\tilde H) =
dim_k \sum_i H^*(X;k) \geq 2^r$.  $\hfill\Box$\\

\begin{Rem} The above results are analogous to results in
homological algebra using a multiplicative structure (if it exists
(!)) on the minimal resolution of a finite module over a
polynomial ring (see [BE],[Av]). As in the latter case it would
suffice for the result of the corollary, to assume that $\tilde H$
carries  an (associative) dg-module structure over $K_r(m)$
(cf.[Av], Prop.6.4.1). Unfortunately we do not know of any new
examples to which the above corollary could be applied, but it
might be interesting to study the question of existence of such
module structures for special classes of spaces, e.g. (rationally)
formal spaces in the sense of Sullivan.

\end{Rem}

I would like to thank Bernhard Hanke  for useful suggestions
concerning a preliminary version of this note.

\vspace{1cm}

\noindent {\footnotesize {\it {Volker Puppe\\
Fachbereich Mathematik und Statistik\\
Universt\"at Konstanz\\
D-78457 Konstanz\\
Germany\\}
 E-mail: volker.puppe@uni-konstanz.de}}

\end{document}